%


\magnification=\magstep1
\def\forces{\parallel\!\!\! -}
\def\restrict{{\restriction}}

\def\Smallskip{\vskip1.4truecm}
\def\Bigskip{\vskip2.2truecm}

\def\qed{{\vcenter{\hrule height.4pt \hbox{\vrule width.4pt height5pt
 \kern5pt \vrule width.4pt} \hrule height.4pt}}}
\def\ok{\vbox{\hrule height 8pt width 8pt depth -7.4pt
    \hbox{\vrule width 0.6pt height 7.4pt \kern 7.4pt \vrule width 0.6pt height 7.4pt}
    \hrule height 0.6pt width 8pt}}
\def\nt{{\leq}\kern-1.5pt \vrule height 6.5pt width.8pt depth-0.5pt \kern 1pt}
\def\sd{{\times}\kern-2pt \vrule height 5pt width.6pt depth0pt \kern1pt}
\def\notin{{\in}\kern-5.5pt / \kern1pt}
\def\ZZ{{\rm Z}\kern-3.8pt {\rm Z} \kern2pt}
\def\NN{{\rm I\kern-1.6pt {\rm N}}}
\def\QQ{{\Bbb Q}}
\def\CC{{\Bbb C}}
\def\DD{{\Bbb D}}
\def\AA{{\Bbb A}}
\def\BB{{\Bbb B}}
\def\EE{{\Bbb E}}
\def\LL{{\Bbb L}}
\def\II{{\Bbb I}}

\def\PP{{\Bbb P}}
\def\RR{{\Bbb R}}
\def\KK{{\rm I\kern-1.6pt {\rm K}}}
\def\11{{\rm 1}\kern-2.2pt {\rm \vrule height6.1pt
    width.3pt depth0pt} \kern5.5pt}
\def\zp#1{{\hochss Y}\kern-3pt$_{#1}$\kern-1pt}

\def\egs{\vrule height 6pt width.5pt depth 2.5pt \kern1pt}
\font\small=cmr8 scaled\magstep0

\font\capit=cmcsc10 scaled\magstep0

\font\sanse=cmss10 scaled\magstep0
\overfullrule=0pt
\openup1.5\jot
\input mssymb

\noindent {\capit Amoeba--absoluteness and projective measurability}
\Bigskip
\noindent J\"org Brendle\footnote{$^*$}{{\small
The author would like to thank the MINERVA-foundation
for supporting him}}
\Smallskip
{\baselineskip=0pt {\small
\noindent 
Abraham Fraenkel Center for Mathematical Logic,
Department of Mathematics,
Bar--Ilan University,
52900 Ramat--Gan, Israel 
\smallskip
and
\smallskip
\noindent Mathematisches Institut der Universit\"at T\"ubingen,
Auf der Morgenstelle 10, 7400 T\"ubingen, Germany}}
\Bigskip
\noindent{\sanse Abstract.} We show that
$\Sigma^1_4$--Amoeba--absoluteness implies that $\forall a \in \RR
\; (\omega_1^{L[a]} < \omega_1^V)$, and hence
$\Sigma^1_3$--measurability. This answers a question of Haim Judah
(private communication).
\Bigskip
{\capit Introduction}
\Smallskip
We study the relationship between Amoeba forcing and projective
measurability. Recall that the {\it Amoeba partial order} $\AA$
is defined as follows.
\smallskip
\centerline{$A \in \AA \Longleftrightarrow A \subseteq 2^\omega
\;\land\; A$ open $\land\; \mu(A) < {1 \over 2}$}
\centerline{$A \leq B \Longleftrightarrow B \subseteq A$}
\smallskip
\noindent Amoeba forcing generically adds a measure one set of
random reals. Its importance in the investigation of measurability
of projective sets stems from the classical result, due to
Solovay, that
\smallskip
\centerline{(*) \hskip 2truecm {\it all $\Sigma^1_2$--sets
are measurable} $\Longleftrightarrow \forall a \in \RR \;
(\mu (Ra (L[a])) = 1)$}
\smallskip
\noindent (see, e.g., [JS 2, 0.1. and $\S$ 3]). Here $Ra(M)$ denotes
the set of reals random over a model $M$ of set theory.
\par
The connection between Amoeba forcing and projective measurability
was made more explicit through Judah's study of absoluteness
between models $V \subseteq W$ of set theory such that $W$ is
a forcing extension of $V$ [Ju].
\smallskip
{\capit Definition} (Judah [Ju, $\S$ 2]). Let $V$ be a universe
of set theory. Given a forcing notion $\PP \in V$ we say that $V$ is 
{\it $\Sigma^1_n - \PP$--absolute} iff for every $\Sigma^1_n$--sentence
$\phi$ with parameters in $V$ we have $V \models \phi$ iff $V^\PP
\models \phi$. (So this is equivalent to saying that $\RR^V \prec_{
\Sigma_n^1} \RR^{V^\PP}$.)
\smallskip
\noindent Note that Shoenfield's Absoluteness Lemma [Je, 
Theorem 98] says that $V$ is
always $\Sigma_2^1-\PP$--absolute. Furthermore, 
using (*), Judah showed [Ju, $\S$ 2]
\smallskip
\centerline{(**) \hskip 2truecm {\it all $\Sigma^1_2$--sets are
measurable in $V$ $\Longleftrightarrow$ $V$ is $\Sigma^1_3 -
\AA$--absolute.}}
\smallskip
Whereas there is no way of getting a characterization of
$\Sigma^1_3$--measurability analogous to (*), (**) suggests
the investigation of the relation between
$\Sigma^1_3$--measurability and $\Sigma^1_4 - \AA$--absoluteness.
The main goal of this note is to establish one implication,
namely that $\Sigma^1_4 - \AA$--absoluteness implies
$\Sigma^1_3$--measurability (Theorem 5 in $\S$ 2). Our tools for
proving this theorem are a partial earlier result of Judah's,
who showed Theorem 5 under the additional assumption that
$\forall a \in \RR \; (\omega_1^{L[a]} < \omega_1^V )$, and 
combinatorial ideas due to Cicho{\'n} and Pawlikowski
[CP], which will eventually yield that Judah's additional
assumption is in fact a consequence of $\Sigma^1_4 -
\AA$--absoluteness ($\S$ 1 and Theorem 4 in $\S$ 2).
\Smallskip
{\sanse Notation.} We shall mostly work with $2^\omega$ or
$\omega^\omega$ instead of $\RR$. ${\cal L}$ denotes
the ideal of Lebesgue measure zero sets, and ${\cal B}$ is the
ideal of meager sets. $\Sigma^1_n ({\cal L})$ stands for {\it
all $\Sigma^1_n$--sets are Lebesgue measurable}; and
$\Sigma^1_n ({\cal B})$ means {\it all $\Sigma^1_n$--sets have
the property of Baire}. For a non--trivial $\sigma$--ideal 
${\cal I} \subseteq P(2^\omega)$, let $add({\cal I})$ be
the size of the smallest family of members in ${\cal I}$ 
whose union is not in ${\cal I}$; $cov({\cal I})$ denotes
the least $\kappa$ such that $2^\omega$ can be covered by
$\kappa$ sets from ${\cal I}$; $unif({\cal I})$ is the
cardinality of he smallest subset of the reals which does
not lie in ${\cal I}$; and $cof({\cal I})$ is the size of
the smallest ${\cal F} \subseteq {\cal I}$ such that
every member of ${\cal I}$ is included in a member of ${\cal F}$.
We always have $add({\cal I}) \leq cov({\cal I}) \leq cof({\cal
I})$ and $add({\cal I}) \leq unif({\cal I}) \leq cof({\cal I})$
(see, e.g., [CP] for details concerning these invariants in
case ${\cal I} = {\cal L}$ or ${\cal B}$).
\par
Our forcing notation is rather standard (see [Je] for any
notion left undefined here). We confuse to some extent
Boolean--valued models $V^\PP$ and forcing extensions $V[G]$,
$G$ $\PP$--generic over $V$. For p.o.s $\PP$, $\QQ$,
$\PP <_c \QQ$ means that $\PP$ can be completely embedded in
$\QQ$. For a sentence of the $\PP$--forcing language $\phi$,
$\Vert \phi \Vert$ is the Boolean value of $\phi$. $\PP$--names
for objects in the forcing extension are denoted by symbols
like $\breve r$. Finally, $\BB$ will stand for the random
algebra, $\CC$ for the Cohen algebra, 
and $\DD$ for the Hechler p.o. (see, e.g., [BJS]).
\Smallskip
{\sanse Acknowledgments.} I am very much indebted to both
Haim Judah (for sharing with me his insight into projective
measurability and motivating me to work in the area) and
Andrzej Ros{\l}anowski (for several stimulating discussions,
concerning mainly the material in $\S$ 1).
\Bigskip
{\capit $\S$ 1. The combinatorial component}
\Smallskip
We start with a straightforward generalization of one version
of the main result
of [CP]. The proof is included for completeness' sake.
\smallskip
{\capit Theorem 1} (Cicho\'n -- Pawlikowski [CP, $\S$ 1]). 
{\it Assume that $\CC \leq_c \PP$, 
and that for any uncountable $T \subseteq \PP$ 
there is an $s \in \CC$ such that for all $\ell \in \omega$
there exists $F \subseteq T$ of size $\ell$ such that 
any $t$ extending $s$ is compatible with $\bigcap F \in \PP$.
Then there is a family $\{ A_x ; \; x \in \omega^\omega \cap V \}$
of Lebesgue measure zero sets in $V^\CC$ 
such that for all $z \in V^\PP$,
$\{ x \in \omega^\omega \cap V ; \; z \not\in A_x \}$ is at most
countable.}
\smallskip
{\it Proof.} Let $\{ \tau_n ; \; n \in \omega \}$ be a one--to--one
enumeration of $\omega^{<\omega}$; set $code(\tau) = n$ iff
$\tau = \tau_n$ for any $\tau \in \omega^{<\omega}$.
Let $\{ C^n (i) ; \; i \in \omega \}$ be an enumeration of
all open intervals in the unit interval $\II
=[0,1]$ with rational endpoints of length
$2^{-n}$, For $x,y \in \omega^\omega$ let 
$$B^n_{x,y} = \cases{C^n (\tau_{y(n)} (code (x \restrict y(n+1))))
&if $code(x \restrict y(n+1)) \in dom(\tau_{y(n)})$ \cr
\emptyset &if not \cr }$$
Let $B_{x,y} = \bigcap_n \bigcup_{m > n} B_{x,y}^m$.
Clearly $\mu (B_{x,y}) = 0$. We claim that if $c$ is Cohen over $V$,
$A_x = B_{x,c}$ for $x \in \omega^\omega \cap V$, then $\{ A_x ; \;
x \in \omega^\omega \cap V \}$ is the required family.
\par
For suppose not. Then there are a $\PP$--name $\breve z$, an
uncountable set $T \subseteq \omega^\omega \cap V$, $T \in V$,
conditions $p_x \in \PP$, and $k_x \in \omega$ ($x \in T$) such
that
\smallskip
\centerline{$p_x \forces_\PP \forall n \geq k_x \; (\breve z \not\in
B_{x,\breve c}^n )$ \hskip 2truecm (*).}
\smallskip
\noindent Choose $T' \subseteq T$ uncountable and 
$k \in \omega$ such that $\forall x \in T' \; (
k_x = k)$. Fix $s \in \CC$ according to $T'$. 
Let $\ell \geq k$, $\ell \geq lh(s)$, and choose
$F \subseteq \omega^\omega$ of size $2^\ell$ such that $ \{ p_x ;
\; x \in F \} $ satisfies the requirements of the Theorem.
Next let $n > \ell$ be
such that $|\{ x \restrict n ; \; x \in F \} | = 2^\ell$.
Let $F = \{ x_i ; \; i \in 2^\ell \}$, and choose $i_0, ..., i_{2^\ell - 1}$
such that $C^\ell (i_0) \cup ... \cup C^\ell (i_{2^\ell - 1}) =
\II$. Let $m \in \omega$ be such that $\tau_m (code (x_0 \restrict n ))
= i_0 , ..., \tau_m (code (x_{2^\ell - 1 } \restrict n )) = i_{2^\ell - 1}$.
Let $t \leq s$ be such that $t(\ell) = m$, $t(\ell + 1) = n$.
Then $\bigcup_{i \in 2^\ell} C^\ell (\tau_{t(\ell)} (code (x_i \restrict
t(l+1)))) = \II$, i.e.
$$t \cap \bigcap \{ p_x; \; x \in F \} \forces_\PP \breve z \in
\bigcup_{i \in 2^\ell} C^\ell (\tau_{\breve c (\ell)} (code (x_i \restrict
\breve c (\ell + 1 )))) = \bigcup_{i \in 2^\ell} B^\ell_{x_i , \breve c},$$
contradicting (*). $\qed$
\bigskip
As each open set in $2^\omega$ can be written as a
countable disjoint union of sets of the form $[\sigma] = \{ f \in 
2^\omega ; \; \sigma \subseteq f \}$, where $\sigma \in 2^{< \omega}$, we
can think of a condition $A$ in the Amoeba
algebra $\AA$ as a function $\phi : \omega \to \bigcup_{
i \in \omega} P(2^i)$ with $\phi (i) \in P( 2^i)$ such that $\sigma \in
\phi (i)$ iff $\sigma \in 2^i$ and $\sigma$ lies in the countable disjoint
decomposition of $A$.
We can furthermore assume that $\phi$ has the property:
\smallskip
(*) \hskip 2truecm $\forall \sigma \in 2^i \setminus \phi (i) \; ( \mu (\cup \{
[\tau] ; \; \tau \supseteq \sigma \; \land \; \exists j > i \; (\tau
\in \phi (j)) \} ) < 2^{-i}).$
\smallskip
\noindent (Then $\phi$ is unique.) 
We define a p.o. $\AA '$ as follows.
$$(u,\phi) \in \AA' \Longleftrightarrow \cases{ {\rm 1)} \; dom (\phi)
= \omega \; \land\; \forall i \in \omega\; (\phi (i) \in P(2^i)) \;\land\;
\phi \;{\rm satisfies \; (*)} \cr
{\rm 2)} \;u \subseteq \phi \;(u \;{\rm  is \; an \; initial 
\; segment \; of } \;\phi) \cr
{\rm 3)} \;\mu (\cup \{ [\sigma] ; \; \exists i \in \omega \; (\sigma 
\in \phi (i) ) \}) < {1 \over 2} \cr}$$
$$(u,\phi) \leq (v,\psi) \Longleftrightarrow u \supseteq v \;\land\;
\forall i \; \forall \sigma \in \psi (i) \; \exists j \leq i \; \exists
\tau \in \phi (j) \; (\sigma \supseteq \tau)$$
\par
{\capit Lemma 1.} {\it $\AA$ and $\AA '$ are equivalent.}
\smallskip
{\it Proof.} We define $\Phi : \AA \to \AA '$ as follows. $\Phi (\phi)
= (u, \phi)$, where $u \subseteq \phi$ is such that $dom (u)$ is
maximal with the following property: for any extension $\psi 
\supseteq \phi$ in $\AA$, $\psi \restrict dom (u) = \phi \restrict 
dom (u)$. We claim that $\Phi$ is a dense embedding. 
\par
Clearly $\psi \leq \phi$ implies $\Phi (\psi) \leq \Phi (\phi)$, and
$\psi \bot \phi$ implies $\Phi (\psi) \bot \Phi (\phi)$. To check
density, choose $(u, \phi) \in \AA '$. Let $i := dom (u) - 1$; and
set $S_\phi := \{ \sigma \in 2^i ; \; $ for no $j \leq i$ does there
exist $\tau \in u(j)$ such that $\sigma \supseteq \tau \}$. For
$\sigma \in S_\phi$ we have $m_\sigma : = \mu ( [\sigma] \setminus
\cup \{ [\tau] ; \; \tau \supseteq \sigma \; \land \; \exists i
\geq dom (u) \; (\tau \in \phi (i) ) \} ) > 0$. Let $a :=
min \{ m_\sigma ; \;\sigma\in S_\phi \}$; and note that $\sum_{\sigma \in
S_\phi} m_\sigma > {1 \over 2}$.
\par
Now define $\psi$ satisfying (*) such that \par
1) $\forall i \in dom (u) \; (\psi (i) = \phi (i))$ \par
2) $\forall i \geq dom (u) \;\forall \tau_1 \in \phi (i) \; \exists
j \leq i \;\exists \tau_2 \in \psi (j) \; (\tau_2 \subseteq \tau_1)$
\par
3) ${1 \over 2} > \mu (\cup \{[\tau]; \; \exists i \in \omega \;
(\tau \in \psi (i)) \} ) > {1 \over 2} - {a \over 2}$
\par
4) for each $\sigma \in S_\phi$, $\mu ([\sigma] \setminus \cup \{
[\tau] ; \; \tau \supseteq \sigma \;\land\; \exists i \geq n \;
(\tau \in \psi (i)) \}) \geq {a \over 2}$
\par
\noindent This is clearly possible. By construction we have $\Phi (\psi)
= (u, \psi) \leq (u , \phi)$. $\qed$
\smallskip
Next define $\AA '' \subseteq \AA '$ by
$$(u,\phi) \in \AA '' \Longleftrightarrow \cases{ {\rm for \; some } 
\; n \in
\omega \; {\rm  we \; have } \;\mu(\cup \{ [\sigma] ; \; \exists i \in dom (u) \;
(\sigma \in u(i)) \}) > {1 \over 2} - {1 \over 2^n}, \cr
\mu (\cup \{ [\sigma] ; \; \exists i \in dom (u)-1 \; (\sigma \in u(i))\})
\leq {1 \over 2} - {1 \over 2^n}, \cr
{\rm and} \;\mu(\cup \{ [\sigma] ; \; \exists i \geq dom (u) \; (\sigma
\in \phi (i) ) \}) < {1 \over 2^{n + 7}}. \cr}$$
Clearly $\AA''$ is dense in $\AA '$. Finally we want to define $h:
\AA '' \to \CC$ giving rise to a complete embedding
of $\CC$ into $\AA$. To this end, let $f: \omega \to \omega$ be such
that $\forall n \; \exists^\infty i \; (f(i) = n)$. For $(u,\phi)
\in \AA ''$ and $n 
\in \omega$ such that ${1\over 2} - {1 \over 2^{n+1}} \geq \mu (\cup
\{ [\sigma] ; \; \exists i \in dom (u) \; (\sigma \in u(i))\}) >
{1 \over 2} - {1 \over 2^n}$ and each $j \leq n$ choose $i_j$ minimal
such that $\mu(\cup\{ [\sigma] ; \; \exists i \in i_j \; (\sigma \in 
u(i))\}) > {1 \over 2} - {1 \over 2^j}$, and let $h((u,\phi)) =
\langle f(i_0) \rangle \hat{\;} ... \hat{\;} \langle f (i_n) \rangle$.
We leave it to the reader to verify that $h : \AA '' \to \CC$ is
indeed a projection (in the forcing theoretic sense). Furthermore,
given $T \subseteq \AA ''$ uncountable we can find $T' \subseteq T$
uncountable and $u$ such that all elements of $T'$ are of the form
$(u,\phi)$ for some $\phi$. Then there is an $s \in \CC$ such that
$\forall (u,\phi) \in T' \; (h((u,\phi)) = s)$. Next, given
$\ell \in \omega$, we can find $F \subseteq T'$ of size $\ell$
such that $\cap F \in \AA ''$. Clearly $h(\cap F) = s$
and so any extension of $s$ in $\CC$ will be compatible with
$\cap F$. Hence we have proved that $\AA ''$ satisfies the requirements
of Theorem 1. Using Lemma 1 we get
\smallskip
{\capit Theorem 2.} {\it There is a family $\{ A_x ; \; x \in \omega^\omega
\cap V \}$ of Lebesgue measure zero sets in $V^\AA$ such that
for all $z \in V^\AA$, $\{ x \in \omega^\omega \cap V ; \; z \not\in
A_x \}$ is at most countable.} $\qed$
\smallskip
{\capit Corollary 1.} {\it Let $V \subseteq W$ be models of $ZFC$ such
that $\omega_1^V = \omega_1^W$. Then there is no real random in $W^\AA
$ over $V^\AA$.}
\smallskip
{\it Proof.} Let $\{A_x ; \; x \in \omega^\omega \cap W \}$ be as in Theorem
2 and note that $\forall z \in \omega^\omega \cap W^\AA \; \exists x \in
\omega^\omega \cap V \; (z \in A_x )$. Hence any real in $W^\AA$ lies
in a measure zero set coded in $V^\AA$. $\qed$
\smallskip
\noindent Using a similar argument as in [CP, $\S$ 3] we can prove
\smallskip
{\capit Corollary 2.} {\it After adding one Amoeba real, $cov ({\cal L}) 
= add ({\cal L}) = \omega_1$ and $unif ({\cal L}) = cof ({\cal L})
= 2^\omega$.} $\qed$
\smallskip
\noindent We note that in [BJS, $\S$ 2] results much stronger than Theorem 2
and the Corollaries were proved for the Hechler p.o. $\DD$; e.g.
it was shown that after adding a Hechler real, $add({\cal B}) =
unif({\cal B}) = \omega_1$ and $cof({\cal B}) = cov({\cal B}) =
2^\omega$ [BJS, 2.5.]. Accordingly we ask:
\smallskip
{\capit Question} [BJS, 2.7.]. {\it Is $unif({\cal B}) = \omega_1$
and $cov({\cal B}) = 2^\omega$ after adding an Amoeba real?}
\bigskip
Before ending this section I wish to include a few comments,
some of which are due to Andrzej Ros{\l}anowski.
\smallskip
{\capit Definition} (implicit in [Tr 2]). A p.o. $\PP$ is
said to have {\it $(\omega_1 , \omega)$--caliber} iff for any
uncountable $T \subseteq \PP$ of size $\omega_1$ there is
a countable $F \subseteq T$ such that $\cap F \in \PP$.
\smallskip
\noindent This is equivalent to: {\it any set of ordinals $A$
in $V^\PP$ of size $\geq \omega_1$ has a countable subset $B$ 
in $V$} [Tr 2]. It is easy to see that if $\CC \leq_c \PP$ and
$\PP$ has $(\omega_1 , \omega)$--caliber, then the assumptions
of Theorem 1 are satisfied. Furthermore the Amoeba algebra
$\AA$ has $(\omega_1 , \omega)$--caliber (the proof for this is
similar to the corresponding proof for the random algebra $\BB$,
given in [Tr 2]). This gives an alternative argument to prove
Theorem 2. --- Our reason for giving the (slightly more difficult)
above argument involving $\AA '$ and $\AA ''$ is that
along the same lines results corresponding to Theorem 2 and the
Corollary can be proved for p.o.s not having
$(\omega_1 , \omega)$--caliber. We include two examples for such
p.o.s:
\smallskip
--- the eventually different reals p.o. $\EE$, due to Miller
[Mi]:
\smallskip
\centerline{$(s,G) \in \EE \Longleftrightarrow s \in
\omega^{<\omega} \;\land\; G \in [\omega^\omega]^{<\omega}$}
\centerline{$(s,G) \leq (t,H) \Longleftrightarrow s \supseteq
t \;\land\; G \supseteq H \;\land\; \forall g \in H \; \forall
i \; (dom(t) \leq i < dom (s) \rightarrow s(i) \neq g(i))$}
\smallskip
--- the localization p.o. $\LL$ (see, e.g., [Tr 3, $\S$ 2]):
\smallskip
\centerline{$(\sigma,G) \in \LL \Longleftrightarrow \sigma \in ([
\omega]^{< \omega})^{< \omega} \;\land\; \forall i \in dom(\sigma)
\; (\vert \sigma (i) \vert = i + 1)\;\land\; G \in
[\omega^\omega]^{\leq dom(\sigma) + 1}$}
\centerline{$(\sigma , G) \leq (\tau,H) \Longleftrightarrow
\sigma \supseteq \tau \;\land\; G \supseteq H \;\land\; \forall 
g \in H \;\forall
i \;(dom(\tau) \leq i < dom(\sigma) \rightarrow g(i) \in \sigma
(i))$}
\smallskip
\noindent Let $\{ f_\alpha ; \; \alpha < \omega_1 \} \subseteq
\omega^\omega$ be a family of pairwise eventually different reals
(i.e. $\alpha \neq \beta \rightarrow \exists n \; \forall k \geq n
\; (f_\alpha (k) \neq f_\beta (k))$). Then $\{ (\langle\rangle ,
\{ f_\alpha \} ) ; \; \alpha < \omega_1 \}$ is an uncountable
set of conditions in $\EE$ (and $\LL$) such that no countable
subset has nontrivial intersection, thus witnessing that
$\EE$ and $\LL$ do not have $(\omega_1 , \omega)$--caliber. We leave
it to the reader to verify that both still satisfy the assumptions
of Theorem 1, however (note that both have a definition similar to,
but easier than, $\AA ''$). \par
(The localization p.o. $\LL$ arose from Bartoszy{\'n}ski's
characterization of the cardinal $add({\cal L})$ [Ba], and is 
closely related to the Amoeba algebra $\AA$. Truss [Tr 3, $\S$ 4]
showed that $\AA <_c \LL$. By the above discussion the converse
cannot hold.)
\Bigskip
{\capit $\S$ 2. The projective part}
\Smallskip
We first introduce a notion closely related to {\it absoluteness},
and discuss the relationship between the two notions.
\smallskip
{\capit Definition} (Judah [Ju, $\S$ 2]). Let $V$ be a universe of set
theory. Given a forcing notion $\PP \in V$ we say that $V$ is {\it $\Sigma^1_n -
\PP$--correct} iff for every $\Sigma^1_n$--formula $\phi (x)$ 
with parameters in $V$ and for
every $\PP$--name $\tau$ for a real we have $V[\tau] \models \phi (\tau)$
iff $V^\PP \models \phi (\tau)$.
\smallskip
{\capit Lemma 2.} {\it Suppose $\PP <_c \QQ$. Then: \par
\item{(i)} $\Sigma^1_n - \QQ$--correctness implies $\Sigma^1_n -
\PP$--correctness. \par
\item{(ii)} $\Sigma^1_{n+1} - \QQ$--absoluteness + $\Sigma^1_n - \QQ
$--correctness implies $\Sigma^1_{n+1} - \PP$--absoluteness.}
\smallskip
{\it Proof.} We prove both (i) and (ii) by induction on $n$. \par
(i) $n = 2$ follows from Shoenfield's Absoluteness Lemma.
Suppose it is true for $n \geq 2$ and assume $V$ is $\Sigma^1_{n+1} - \QQ
$--correct. Let $\phi (x)$ be a $\Sigma^1_{n+1}$--formula, $\phi (x) =
\exists y \psi (y,x)$ where $\psi$ is $\Pi_n^1$. Suppose first that
$V[\tau] \models \phi (\tau)$. Then $V[\tau] \models \exists x \psi (x ,
\tau )$. So there is a $\PP$-name $\sigma$ such that $V[\tau] = V[\sigma,
\tau] \models \psi (\sigma , \tau)$. By induction $V^\PP \models \psi
(\sigma,
\tau)$; thus $V^\PP \models \phi (\tau)$. \par
Assume now that $V^\PP \models \phi (\tau)$. Hence $V^\PP \models \exists
x \psi (x, \tau)$; and we can again find a $\PP$--name $\sigma$ such
that $V^\PP \models \psi (\sigma, \tau)$. By induction $V[\sigma , \tau ]
\models \psi (\sigma , \tau )$. So $\Sigma^1_n - \QQ$--correctness implies
$V^\QQ \models \psi (\sigma , \tau)$; thus $V^\QQ \models \phi (\tau)$.
Hence by $\Sigma^1_{n+1} - \QQ$--correctness $V[\tau] \models \phi
(\tau)$. \par
(ii) $n = 1$ follows from Shoenfield's Absoluteness Lemma. Suppose (ii)
is true for $n \geq 1$ and assume $V$ is $\Sigma^1_{n+2} - \QQ$--absolute
and $\Sigma^1_{n+1} - \QQ$--correct. By (i) $V$ is also $\Sigma^1_{n+1} -
\PP$--correct. Let $\phi$ be a $\Sigma^1_{n+2}$--sentence, $\phi
= \exists x \psi (x)$, where $\psi$ is $\Pi^1_{n+1}$. Suppose first that 
$V \models \phi$; i.e. $V \models \psi (a)$ for some $a \in V$. By induction
$V^\PP \models \psi (a)$; thus $V^\PP \models \phi$. \par
Assume now that $V^\PP \models \phi$; i.e. $V^\PP \models \psi (\tau)$
for some $\PP$--name $\tau$. By $\Sigma^1_{n+1} - \PP$--correctness 
$V[ \tau] \models \psi (\tau)$. Hence $\Sigma^1_{n+1} - \QQ$--correctness
implies $V^\QQ \models \phi$. Thus $V \models \phi$ by $\Sigma^1_{n+2} -
\QQ$--absoluteness. $\qed$
\smallskip
{\capit Lemma 3} (Truss [Tr 1, 6.5]). {\it $\DD <_c \AA$.} $\qed$
\smallskip
{\capit Definition} (Judah -- Shelah [JS 1, $\S$ 0]). 
A ccc notion of forcing $(\PP,\leq)$ is called {\it
Souslin} iff it can be thought of as a $\Sigma^1_1$--subset
of the reals $\RR$ with both $\leq$ and $\bot$ (incompatibility)
being $\Sigma^1_1$--relations (in the plane $\RR^2$).
\smallskip
\noindent Note that all p.o.s discussed in this paper are
Souslin.
\smallskip
{\capit Theorem 3} (Judah [Ju, $\S$ 2]). {\it Assume that $\forall a \in \RR
\; (\omega_1^{L[a]} < \omega_1^V)$, and $\PP \in V$ is a Souslin forcing.
Then $V$ is $\Sigma^1_3 - \PP$--correct.} $\qed$
\bigskip
{\capit Theorem 4.} {\it $\Sigma^1_4 - \AA$--absoluteness implies that 
$\forall a \in \RR$ $(\omega_1^{L[a]} < \omega_1^V)$.}
\smallskip
{\capit Corollary 3.} {\it $\Sigma^1_4 - \AA$--absoluteness implies $\Sigma^1_3
- \AA$--correctness, $\Sigma^1_4 - \DD$--absoluteness, and $\Sigma^1_3 -
\DD$--correctness.}
\smallskip
{\capit Theorem 5.} {\it $\Sigma^1_4 - \AA$--absoluteness implies $\Sigma^1_3
({\cal L}) $ and $\Sigma^1_3 ({\cal B})$.}
\smallskip
The proof of Theorem 4 follows the lines of the proof of 2.6 in [BJS].
Theorem 5 is a consequence of Theorem 4 and a result in [Ju, $\S$ 2].
We give the proof here for completeness' sake. --- Note that 
$\Sigma^1_3 - \DD$--absoluteness is equivalent to $\Sigma^1_2 (
{\cal B})$ [Ju, $\S$ 2]. Thus the implication {\it $\Sigma^1_3 -
\AA$--absoluteness $\Longrightarrow \;\Sigma^1_3 - \DD$--absoluteness}
(immediate from Lemmata 2 and 3) is just another version of the
Raisonnier--Stern Theorem; and Corollary 3 may be thought of as the
corresponding result for $\Sigma^1_4$.
\smallskip
{\it Proof of Theorem 4.} Suppose there is an $a \in \RR$ such that
$\omega_1^{L[a]} = \omega_1^V$. By $\Sigma^1_3 - \AA$--absoluteness we have
$\Sigma^1_2 ({\cal L})$; i.e. $\forall b \in \RR \; (\mu (Ra (L[b])) =1)$
(see the beginning of this section). Note that $x \in Ra (L[b])$ is
equivalent to
\smallskip
\centerline{$\forall c \; (c \not\in L[b] \cap BC \; \lor \; \hat c$ is not
null $\lor\; x \not\in \hat c)$,}
\smallskip
\noindent where $BC$ is the set of Borel codes which is $\Pi^1_1$
[Je, Lemma 42.1], 
and for $c \in BC$, $\hat c$ is the set coded by $c$. As $L[b]$
is $\Sigma^1_2$ [Je, Lemma
41.1], $Ra(L[b])$ is a $\Pi^1_2$--set. Hence $\forall
b \in \RR \; (\mu (Ra (L[b])) = 1)$ which is equivalent to
\smallskip
\centerline{$\forall b \exists c \; (c \in BC \; \land \; \hat c$ is null
$\land \;\forall x \; (x \in \hat c \; \lor \; x \in Ra(L[b])))$}
\smallskip
\noindent is a $\Pi^1_4$--sentence. So it is true in $V^\AA$ by
$\Sigma^1_4$--absoluteness; in particular $Ra(L[a][r])$ (where $r$
is Amoeba over $V$) has measure
one in $V[r]$ which implies that there is a random real in $V[r]$
over $L[a][r]$, contradicting Corollary 1 in $\S$ 1. $\qed$
\smallskip
{\it Proof of Corollary 3.} Follows from Theorems 3 and 4 and Lemmata
2 and 3. $\qed$
\smallskip
{\it Proof of Theorem 5} (Judah). Let $\phi (x)$ be a $\Sigma^1_3$--formula and
$A = \{ x; \; \phi (x) \}$. We shall show that $A$ is measurable in $V$.
First note that the sentence {\it $A$ has measure zero} is equivalent to
\smallskip
\centerline{$\exists c \; (c \in BC \; \land \; \mu (\hat c) = 0 \;
\land \; \forall x \; ( \neg \phi (x) \; \lor \; x \in \hat c ))$,}
\smallskip
\noindent which is $\Sigma^1_4$. So by $\Sigma^1_4
- \AA$--absoluteness, if $A$ is null
in $V^\AA$, it is also null in $V$. \par Hence assume that $A$ is not null
in $V^\AA$. As $\mu (Ra (V)) = 1$ in $V^\AA$, there is $r \in Ra (V) \cap
A$ in $V^\AA$; i.e. $V^\AA \models \phi (r)$. By $\Sigma^1_3 -
\AA$--correctness $V[r] \models \phi (r)$. Now let $\phi (x) = \exists y
\psi (x,y)$, where $\psi$ is $\Pi^1_2$. Then there is an $s \in V[r]$
such that $V[r] \models \psi (r,s)$. If $a \in V$ codes the parameters of
$\psi$ and of the name of $s$, we have by Shoenfield's Absoluteness Lemma
$L[a][r] \models \psi (r,s)$. Let $\breve r$ be the $\BB$--name for the
random real $r$ and $s(\breve r)$ a $\BB$--name for $s$. Then
the Boolean value $\Vert \psi (\breve r , s(\breve r)) \Vert$ is non--zero.
Furthermore, if $r' \in \Vert \psi (\breve r , s(\breve r )) \Vert
\cap V$
is random over $L[a]$, then $L[a] [r'] \models \psi (r', s (r'))$ and ---
by absoluteness --- $V \models \psi(r', s(r'))$;
in particular $V \models \phi (r')$. \par
By $\Sigma^1_3 - \AA$--absoluteness we have that $\mu (Ra (L [a])) = 1$
in $V$ (cf Introduction). And the previous discussion gives
us that $Ra(L[a]) \cap \Vert \psi(\breve r, s(\breve r)) \Vert \subseteq A$.
This shows that any non--null $\Sigma^1_3$--set has positive inner
measure; and it is easy to conclude from this that any $\Sigma^1_3$--set
is indeed measurable. 
\par
Finally, $\Sigma^1_3 ({\cal B})$ follows along the same lines because
$\AA$ adds a comeager set of Cohen reals. $\qed$
\smallskip
{\capit Questions.} {\it 1) Does $\Sigma^1_3 ({\cal L})$ imply
$\Sigma^1_4 - \AA$--absoluteness? \par
2) Does $\Sigma^1_4$--Amoeba--meager--absoluteness (or $\Sigma^1_4
- \DD$--absoluteness) imply $\Sigma^1_3 ({\cal B})$?} (cf [Tr 1,
$\S$ 5] for Amoeba--meager forcing ---
the problem here is whether
$\Sigma^1_4$--Amoeba--meager--absoluteness implies
$\forall a \in \RR \; (\omega_1^{L[a]} < \omega_1^V)$; 
cf [BJS, $\S$ 2] for $\DD$ --- the problem here is that $\DD$
does not add a comeager set of Cohen reals) \par
{\it 3) Does $\forall n \; (V$ is $\Sigma^1_n - \AA$--absolute $)$
imply projective measurability?} \par
{\it 4) {\rm (Judah)} 
Does $\Sigma^1_3 ({\cal L})$ imply $\Sigma^1_3 ({\cal B})$?
} (cf Corollary 3)
\Bigskip
{\capit References}
\Smallskip
\itemitem{[Ba]} {\capit T. Bartoszy\'nski,} {\it Combinatorial
aspects of measure and category,} Fundamenta Mathematicae, vol. 127
(1987), pp. 225-239.
\smallskip
\itemitem{[BJS]} {\capit J. Brendle, H. Judah and S. Shelah,}
{\it Combinatorial properties of Hechler forcing,}
preprint.
\smallskip
\itemitem{[CP]} {\capit J. Cicho\'n and J. Pawlikowski,} {\it
On ideals of subsets of the plane and on Cohen reals,}
Journal of Symbolic Logic, vol. 51 (1986), pp. 560-569.
\smallskip
\itemitem{[Je]} {\capit T. Jech,} {\it Set theory,} Academic Press,
San Diego, 1978.
\smallskip
\itemitem{[Ju]} {\capit H. Judah,} {\it Absoluteness for projective
sets,} to appear in Logic Colloquium 1990.
\smallskip
\itemitem{[JS 1]} {\capit H. Judah and S. Shelah,} {\it 
Souslin forcing,} Journal of Symbolic Logic, vol. 53 (1988),
pp. 1188-1207.
\smallskip
\itemitem{[JS 2]} {\capit H. Judah and S. Shelah,} {\it 
$\Delta_2^1$--sets of reals,} Annals of Pure and Applied Logic,
vol. 42 (1989), pp. 207-223.
\smallskip
\itemitem{[Mi]} {\capit A. Miller,} {\it Some properties of measure
and category,} Transactions of the American Mathematical Society,
vol. 266 (1981), pp. 93-114.
\smallskip
\itemitem{[Tr 1]} {\capit J. Truss,} {\it Sets having calibre 
$\aleph_1$,} Logic Colloquium 76, North-Holland, Amsterdam,
1977, pp. 595-612.
\smallskip
\itemitem{[Tr 2]} {\capit J. Truss,} {\it The noncommutativity
of random and generic extensions,} Journal of Symbolic Logic,
vol. 48 (1983), pp. 1008-1012.
\smallskip
\itemitem{[Tr 3]} {\capit J. Truss,} {\it Connections between
different Amoeba algebras,} Fundamenta
Mathematicae, vol. 130 (1988), pp. 137-155.

\vfill\eject\end